\documentclass[10pt,a4paper,twoside]{article}
\usepackage{a4,amsfonts}

\newcommand{\Z}{{\mathbf{Z}}}   
\newcommand{\Q}{{\mathbf{Q}}}
\def\qed{\hspace{\fill}\hbox{${\vcenter{\vbox{              
	\hrule height 0.4pt\hbox{\vrule width 0.4pt height 6pt
	\kern5pt\vrule width 0.4pt}\hrule height 0.4pt}}}$}}

\newtheorem{theorem}{Theorem}
\newtheorem{lemma}{Lemma}[section]
\newtheorem{conjecture}[lemma]{Conjecture}
\newtheorem{definitionhelp}[lemma]{Definition}
\newtheorem{remarkhelp}[lemma]{Remark}
\newtheorem{questionhelp}[lemma]{Question}
\newtheorem{corollary}[lemma]{Corollary}

\newenvironment{definition}{\begin{definitionhelp}\rm}{\end{definitionhelp}}
\newenvironment{remark}{\begin{remarkhelp}\rm}{\end{remarkhelp}}

\def\runninghead#1#2{\pagestyle{myheadings}
\markboth{{\protect\footnotesize\it{\quad #1}}\hfill}
{\hfill{\protect\footnotesize\it{#2\quad}}}}
\input epsf
\newcommand{\epsfig}[1]{\noindent
	\epsfbox{#1}}
\newcommand{\cepsfig}[1]{\par\smallskip\noindent
	\centerline{\epsfig{#1}}\smallskip}

\newcommand{\eqnfig}[2]{
	$$\begin{array}{c}
	\epsfig{#1}
	\end{array} 
	\eqno{\rm (#2)}
	$$}

\begin{document}

\runninghead{\quad Jan A.~Kneissler}{Relations in the Algebra $\Lambda$ \quad}
\title{{\large\it On Spaces of connected Graphs II}\\[0.2cm] \smallskip Relations in the algebra $\Lambda$}
\author{Jan A.~Kneissler}
\date{\quad} 
\maketitle
\thispagestyle{empty}

\vspace{6pt}

\begin{abstract}
{
\noindent
The graded algebra $\Lambda$ defined by Pierre Vogel is of general interest 
in the theory of finite-type invariants of knots and of $3$-manifolds
because it acts on spaces of connected graphs subject to relations 
called IHX and AS.
We examine a subalgebra $\Lambda_0$ 
that is generated by certain elements called $t$ and $x_n$ with $n \geq 3$.
Two families of relations in $\Lambda_0$ are derived and it is shown that the dimension
of $\Lambda_0$ grows at most quadratically with respect to degree. 
Under the assumption that $t$ is not a zero divisor in $\Lambda_0$, 
a basis of $\Lambda_0$ and an isomorphism from $\Lambda_0$ to a 
sub-ring of $\Z[t,u,v]$ is given.
}
\end{abstract}

\section{Main Results}
\label{sectmain}

The theory of Vassiliev invariants, initiated by V.~A.~Vassiliev in \cite{Va},
leads to the study of modules of Chinese characters, i.e.
embedded graphs with univalent and oriented trivalent
vertices, subject to the following relations (the
reader not familiar with these concepts is referred to \cite{BN}, for instance).
\eqnfig{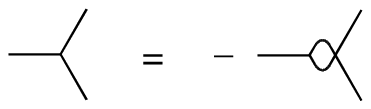}{AS}
\eqnfig{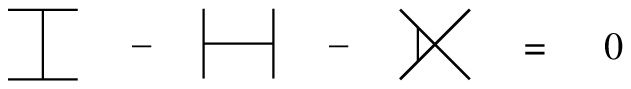}{IHX}

\bigskip
\noindent
\begin{definition} 
Let $F(n)$ denote the $\Q$-vector space of connected tri-/univalent
graphs having exactly $n$ univalent vertices that are labeled with the numbers $1$ to $n$, 
modulo the relations (AS) and (IHX).
If exchanging any two of the labels of $a \in F(n)$ turns $a$ into $-a$,
then $a$ is called {\sl antisymmetric}.
The subspace of $F(3)$ spanned by antisymmetric elements is named $\Lambda$.
\end{definition}

\bigskip
\noindent
All our statements are valid over $\Z[\frac16]$, but for convenience 
we work over $\Q$.
First let us summarize some facts about $\Lambda$ (see \cite{Vo}
for details):
It acts on modules of connected tri-/univalent graphs
by insertion at a trivalent vertex;
(AS) and (IHX) guarantee that 
the result does not depend on the choice of the trivalent vertex. 
So it also acts on itself, which
turns $\Lambda$ into a graded commutative algebra. The following family of graphs 
named $t,x_3,x_4,x_5,\ldots$, represent elements of $\Lambda$ ($x_n$ is homogeneous
of degree $n$ and $t$ has degree $1$).
\cepsfig{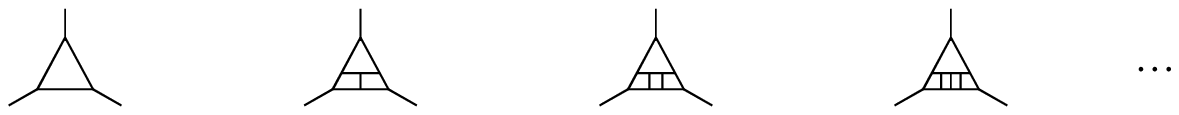}

\noindent
We set $x_1 := 2t$, $x_2 := t^2$ and let $\Lambda_0$ denote the subalgebra of $\Lambda$ that is generated these
elements $x_n$ with $n \geq 1$.
To our present knowledge there are no
counterexamples to the following two conjectures: 

\begin{conjecture}
\label{conj1}
The algebra $\Lambda$ is generated by the $x_n$, i.e.~$\Lambda_0 = \Lambda$.
\end{conjecture}

\begin{conjecture}
\label{conj2}
The element $t$ is not a zero divisor, i.e. the map $\Lambda_0 \rightarrow \Lambda_0$, $\lambda \rightarrow t\,\lambda$
is injective.
\end{conjecture}

\noindent
Considering the results presented here, the verification 
of these conjectures might be the only
remaining hurdle on the way to a complete understanding of $\Lambda$.

\begin{remark}
Computations that are based on results given in the subsequent
paper \cite{cc3} show that Conjecture \ref{conj1} is true 
in degrees $\leq 12$.
Due to Remark \ref{remlbound}, Conjecture \ref{conj2} is true up to degree
$14$.
\end{remark}

\newcommand{\alp}[2]{x_{i+#1}^{}\,x_{j+#2}^{}}
\newcommand{\al}[1]{x_i^{}\,x_{j+#1}^{}}
\newcommand{\bet}[1]{x_{j+#1}^{}}
\medskip
\noindent
To state the main result of this note, we define two families $P_{i,j}$, $Q_{i,j}^k$
of polynomials in the variables $t$ and $x_n$:
\begin{eqnarray*}
P_{i,j}^{} & := & 3\,\alp24 - 3\,t\,\alp14 - 9\,t\,\alp23 - 6\,t^2\,\al4
 + 9\,t^2\,\alp13 + 18\,t^{3}\,\al3  \\
&& - 2\, (x_3 + 2\,t^3)\,\alp12
+ 4\,t\,(x_3 - 4\,t^3)\,\al2 + 8\,t^2\,(t^3 - x_3)\,\al1
 + 3\,t^{i+2}\,\bet4 \\
&& - 9\,t^{i+3}\,\bet3 + (7\,t^3 - x_3)\,t^{i+1}\,\bet2
 + 3\,(x_3 - t^3)\,t^{i+2}\,\bet1 + 2\,(t^3 - x_3)\,t^{i+3}\,x_j^{}
\end{eqnarray*}

\noindent
$Q_{i,j}^k$ is defined recursively for $i,j,k \geq 0$ by

\noindent
\begin{eqnarray*}
Q_{-1,2}^k\!\!  & := & t^k, \\
Q_{0,j}^k\,\,   & := & x_{j+k}^{}, \\
Q_{i+1,j}^k\!\!\!\!\, & := & t\,Q_{i,j+1}^k - \frac12\,Q_{i,j+2}^k
                 + \frac12\,t^j\,Q_{i,2}^k
                 + \frac12\,(x_{j+2}^{}-t^{j+2})\,Q_{i-1,2}^k.
\end{eqnarray*}

\noindent
\begin{theorem}
\label{theo}
The equations $P_{i,j} = P_{j,i}$ for $i,j \geq 1$, $\;Q_{i,1}^1 = Q_{i-1,2}^2$ for $i \geq 0$
are identities in $\Lambda$.
\end{theorem}

\begin{remark}
These relations in $\Lambda$ are homogeneous of degree $i+j+6$ and $2i+2$, respectively.
The relation $Q_{i,1}^1 = Q_{i-1,2}^2$ allows us to
express $x_{2i+2}$ as a polynomial in $t,x_3,x_4\ldots,x_{2i+1}$.
It has already been shown in \cite{Vo} that $\Lambda_0$ is generated by the $x_n$
with odd indices $n$, so the existence of the $Q$-relations is no surprise.

The first (with respect to degree) $P$-relation that is not a consequence of the $Q$-relations
is $P_{1,3} = P_{3,1}$; if all $x_{2i}$ in it are eliminated by $Q$-relations,
one obtains exactly the relation $P$ that has been established in \cite{Vo} and
proven in \cite{Kn}.
\end{remark}

\noindent
It is not known whether the relations of Theorem \ref{theo} together with $x_1 = 2t$ generate the
kernel of the map from the polynomial algebra $\Q[t,x_1,x_2,\ldots]$
to $\Lambda$, but it will be shown in section \ref{sectcor} that any further
relation would be a multiple of $t$, in contradiction to Conjecture \ref{conj2}.

\begin{corollary}
\label{corspan}
The algebra $\Lambda_0$ is spanned (as $\Q$-vector space) by the following set of
monomials $$M := \{ \; t^i,\;t^i\,x_{2n+1}^j,\;t^i\,x_{2n+1}^j\,x_{2n+3}^k \;\,
\vert\;\, i\geq0,\;\; j,\,k,\,n > 0\;\}.$$
If Conjecture \ref{conj2} is true then $M$ is a basis of $\Lambda_0$.
\end{corollary}

\begin{corollary}
\label{coro2}
If $\lambda_d$ denotes the dimension of the part in degree $d$ of $\Lambda_0$, then
$$\left\lfloor \frac{7d}{6} \right\rfloor-3 \; \leq\; \lambda_d \;\leq\; \left\lfloor \frac{d^2}{12} + \frac12
\right\rfloor + 1.$$
The upper bound is saturated if and only if Conjecture \ref{conj2} is valid.
\end{corollary}

\begin{remark}\quad\par
\label{remlbound}
\vspace{-4pt}
\begin{enumerate}
\item[a)]
Using the eight characters of \cite{Vo} and Corollary \ref{corspan}, we verified
that Conjecture \ref{conj2} is true for all $\lambda$ of degree $\leq 14$. This implies that 
the given upper bound is equal to $\lambda_d$ for $d \leq 15$.
\item[b)]
Taking all characters into account, Pierre Vogel improved the lower bound:
$\lambda_d \geq a_d$ with $\sum a_dx^d = \frac1{1-x} +\frac{x^3-x^{16}-x^{23}+x^{26}}{(1-x)(1-x^2)(1-x^3)}$.
This leads to $\lambda_d \geq \lfloor\frac{5d}3\rfloor-2$ for $d \geq 21$.
\end{enumerate}
\end{remark}

\begin{theorem}
\label{theorek}
Conjecture \ref{conj2} is equivalent to the following statement:
There exists a unique ``universal" character $\chi: \Lambda_0 \rightarrow \Q[t,u,v]$
satisfying
$\chi(t) = t$, $\chi(x_0) = 0$, $\chi(x_1) = 2t$, $\chi(x_2) = t^2$ and
for $n\geq 0$:
$$\chi(x_{n+3}) = t\,\chi(x_{n+2}) + u\,\big(2\,\chi(x_{n+1}) - t^{n+1}\big)
+ v\,\big(2\,\chi(x_{n}) - t^n - 2\,(2t)^n\big).$$ 
If it exists, then $\chi$ is an isomorphism of graded algebras between $\Lambda_0$ and the sub-ring
of $\Q[t,u,v]$ that consists of all elements of the form $a + (t^3-tu+v)b$
with $a \in \Q[t]$ and $b \in \Q[t,u,v]$.
\end{theorem}

\begin{remark}\quad\par
\vspace{-4pt}
\begin{enumerate}
\item[a)]
The recursion already appeared in \cite{Vo}, where it is shown 
that for certain values of $t,u,v$ it leads to well-defined homomorphisms 
of $\Lambda_0$, corresponding to the characters coming from 
Lie (super)algebras.
Thus all these characters may be regarded as specializations of the 
(yet hypothetical) universal character.
\item[b)]
We have been informed by Pierre Vogel that $\chi$ might be
derived from a speculative object called ``the universal Lie algebra" 
(see \cite{Vo2}); this may open an approach to Conjecture \ref{conj2}.
\end{enumerate}
\end{remark}

\noindent
{\large\bf Acknowledgments}

\medskip
\noindent
After the computer calculations described in \cite{Kn} had revealed that $\Lambda_0$
is not a polynomial algebra, it was natural to ask if the identity found
in degree $10$ is member of a whole family of relations in $\Lambda$.
In this sense, the main stimulus had been provided by Pierre Vogel; he noticed that 
an identity of values of characters $-$ used in the proof of 
Lemma 7.5 in his paper \cite{Vo} $-$ can be modified to become an element 
of $F(4)$. If it is trivial, then one gets a one-parameter family of
relations in $\Lambda_0$.
Trying to prove this $F(4)$-relation, 
we realized that there should exist a $F(6)$-generalization of it,
which produces a two-parameter family of relations in $\Lambda_0$.
This idea finally lead to the investigations that are presented in \cite{cc1}.

I am also thankful to Jens Lieberum for advises concerning the manuscript,
and to the Studienstiftung for financial support.
A special thank is addressed to all participants of the ``Knot Theory Week"
in Bonn in July '97  for creating an animating atmosphere which promoted the arising
of the ideas that led to these results.

\section{The $P$-relations}

All $P$-Relations are consequences of a single relation in $F(6)$ that
has been proven in \cite{cc1}. To describe it, let $\alpha_{rs}$ denote elements
of $F(6)$ that are given by the following picture:
\cepsfig{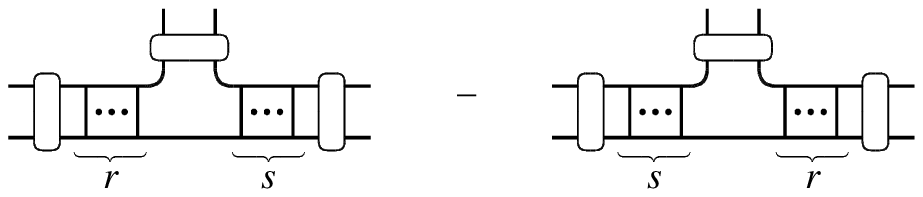}

\noindent
The numbers of vertical edges in the dotted parts are assumed to coincide with
$r$  and $s$ as indicated. The oval boxes are supposed to be expanded in the 
following manner:
\cepsfig{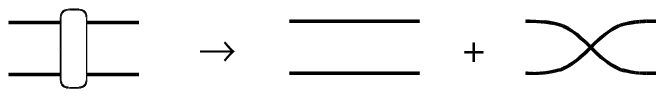}

\noindent
Thus $\alpha_{rs}$ is a linear combination of $16$ graphs, each of them has $2(r+s)$ trivalent
vertices. The univalent vertices are numbered clockwise, starting with the
one on the lower left side. Then $\alpha_{r+1,s+1}$ is the same linear combination of graphs
as the elements $[y]_r^s$ introduced in \cite{cc1}. 
So we may rewrite  equation (39) of \cite{cc1} in the new terminology:
\newcommand{\yp}[2]{\alpha_{#1#2}}
\begin{eqnarray}
\nonumber
3\yp35 &=& 9t\yp34 +3t\yp25 -9t^2\yp24 +6t^2\yp15 +(4t^3+2x_3)\yp23 \\
\label{alph}
&&-18t^3\yp14 +4t(4t^3-x_3)\yp13 +8t^2(x_3-t^3)\yp12
\end{eqnarray}
Now we consider the following operation for graphs of $F(6)$: 
glue $i$ edges between univalent vertices $1$ and $2$ and $j$ edges between $5$ and $6$;
then identify the univalent vertices $2$, $3$ and $5$ to form a new trivalent
vertex:
\cepsfig{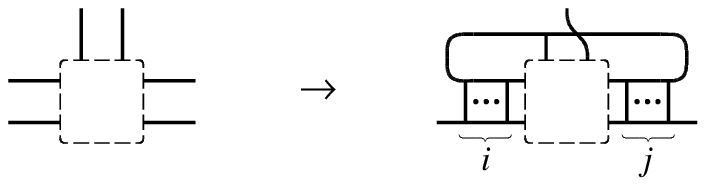}

\noindent
This operation extends to a linear map $\Gamma_{i,j} : F(6) \rightarrow F(3)$. 
It is easy to verify that for $r,s \geq 0$, $i,j \geq 1$ we have 
$$ \Gamma_{i,j}(\alpha_{rs}) \;\;=\;\; 2(2x_{i+r}-t^{i+r})(t^{j+s}-2x_{j+s}) 
- 2(2x_{i+s}-t^{i+s})(t^{j+r}-2x_{j+r}).$$
Finally, by applying $\Gamma_{i,j}$ to (\ref{alph}), we obtain the desired relation $P_{i,j} = P_{j,i}$ in $\Lambda_0$. 

\section{The $Q$-relations}
\newcommand{\W}[1]{\langle#1\rangle}
We use the calculus that has been established in section 5 of \cite{Vo}.
Elements in $\Lambda_0$ are presented by words of length $\geq 2$ in three letters $1$, $2$, $3$.
To the reader's convenience, we will repeat the relations of \cite{Vo}
that will be used:
\vspace{-2ex}
\begin{displaymath}
\end{displaymath}
\begin{displaymath}
\W{12} = -t,\quad \W{12^n1} = x_n
\end{displaymath}
\begin{displaymath}
\forall\sigma \in S_3:\;\; \W{w} = \W{\sigma(w)}
\end{displaymath}
\begin{displaymath}
\W{1v} +\W{2v} +\W{3v} = \W{v1} +\W{v2} +\W{v3} = 0
\end{displaymath}
\begin{displaymath}
\W{u1v} +\W{u2v} +\W{u3v} = 2t\,\W{uv}
\end{displaymath}
\begin{displaymath}
\W{\gamma^2 v} = t\,\W{\gamma v},\quad \W{v \gamma^2} = t\,\W{v \gamma}
\end{displaymath}
\begin{displaymath}
\W{u\gamma v} + \W{u\gamma \sigma_\gamma(v)} = \W{u \gamma} \, \W{\gamma v}
\end{displaymath}
Here $w$ represents a word of length $\geq 2$, $u,v$ words of length $\geq 1$,
$\gamma$ a word of length $1$ and $\sigma_\gamma$ is the transposition of $S_3$ that keeps
$\gamma$ fixed; $\gamma^n$ denotes a word consisting of $n$ copies of the letter $\gamma$.

\bigskip
\noindent
First let us deduce some simple relations:
\begin{eqnarray}
\label{eqg}
\W{u212} &= &-\W{u211} -\W{u213} \;\;=\;\; -\W{u211} -2t\,\W{u23} + \W{u223} + \W{u233}
\nonumber
\\
& = & -t\,\W{u21} -t\,\W{u23} - \W{u222} - \W{u221} \;\;=\;\; -\W{u221}
\\
\nonumber
\W{212^n1} & = & -\W{112^n1} - \W{312^n1} \;\;=\;\;
-t\,\W{12^n1}-2t\,\W{32^n1} + \W{32^{n+1}1} + \W{332^n1} \\
\label{eqe}
&=& -t\,\big(\W{32^n1}+\W{12^n1}+\W{22^n1}\big) - \W{12^{n+1}1}
 \;\;=\;\; - x_{n+2}^{}
\\
\nonumber
\W{232^n1} & = & -\W{332^n1} -\W{132^n1} \;\;=\;\; -t\,\W{32^n1} -2t\,\W{12^n1}
+\W{112^n1} + \W{12^{n+1}1}
\\
\label{eqd}
& = & t\,\W{22^n1} + \W{12^{n+1}1} \;\;=\;\; x_{n+2}^{} - t^{n+2}
\\
\label{eqa}
\W{u212^n1} & = & 2t\,\W{u2^{n+1}1} - \W{u2^{n+2}1} - \W{u232^n1}
\\
\label{eqb}
\W{u232^n1} & = & \W{u2}\,\W{232^n1} - \W{u212^n3}
\;\;=\;\; \W{u2}\,\W{232^n1} + \W{u212^n1} + \W{u212^{n+1}}
\end{eqnarray}
The difference of the last two equations (\ref{eqa}) - (\ref{eqb}) leads to
\begin{eqnarray}
\label{eqc}
\W{u212^n1} & = & t\,\W{u2^{n+1}1} - \frac12\W{u2^{n+2}1} -\frac12t^n\W{u212}
-\frac12\W{u2}\,\W{232^n1}
\end{eqnarray}
For $i,j \geq 0$, $k \in \{1,2\}$ let $q_{i,j}^k := \W{2^{k-1}(12)^i12^j1}$,  then (\ref{eqc}),
(\ref{eqg}), (\ref{eqd}) imply for $i\geq 2$, $j \geq 1$:
\begin{eqnarray*}
q_{i,j}^k &=&t\,q_{i-1,j+1}^k - \frac12q_{i-1,j+2}^k - \frac12\,t^j\,\W{2^{k-1}(12)^{i+1}}
- \frac12\W{2^{k-1}(12)^i}\,\W{232^j1} \\
&=& t\,q_{i-1,j+1}^k - \frac12q_{i-1,j+2}^k + \frac12\,t^j\,q_{i-1,2}^k
+ \frac12 q_{i-2,2}^k \, (x_{j+2}^{}-t^{j+2})
\end{eqnarray*}
This is exactly the recursive formula in the definition of $Q_{i,j}^k$.
We have $q_{0,j}^1 = \W{12^j1} = x_{j+1} = Q_{0,j}^1$ and (\ref{eqe})
$q_{0,j}^2 = \W{212^j1} = -x_{j+2} = -Q_{0,j}^2$
for $j \geq 1$. A little calculation
shows that $q_{1,j}^1 = Q_{1,j}^1$ and $q_{1,j}^2 = -Q_{1,j}^2$ as well.
Thus the equality $Q_{i,j}^k = (-1)^{k+1}q_{i,j}^k$ holds in $\Lambda_0$ for
all $i \geq 0$, $j \geq 1$ and $k \in \{0,1\}$.

\noindent
Because of equation (\ref{eqg}) we have
$$
q_{i-1,2}^2 \;\;=\;\; \W{2(12)^{i-1}1221} \;\;=\;\; -\W{2(12)^{i-1}1212}
\;\;=\;\; -\W{(21)^i212} \;\;=\;\; -\W{(12)^i121}\;\;=\;\;-q_{i,1}^1.
$$
This proves that for $i\geq 1$ the relation $Q_{i,1}^1 = Q_{i-1,2}^2$ is
valid in $\Lambda_0$. It remains to check the
case $i = 0$ which is just the identity $x_2 = t^2$.

\section{Implications}
\label{sectcor}

\begin{lemma}
\label{lembij}
There exists a degree-preserving bijection $\beta$ between the following sets of monomials $\{\; x_{2n+1}^i\,x_{2n+3}^j \,\,
\vert\,\, i,n > 0,\; j\geq 0\;\}$ and $\{\; u^k\,v^l \,\,
\vert\,\, k \geq 0,\; l > 0\;\}$, where the degree of $x_m$, $u$, $v$ are assumed
to be $m$, 2, 3, respectively.
\end{lemma}
\noindent
{\bf Proof }
The map given by $\beta(x_{2n+1}) := u^{n-1}v$ and $\beta(p\,q) := \beta(p)\,\beta(q)$
obviously relates monomials of the same degree. To show that $\beta$ is a bijection,
one has to verify that for any integers $k \geq 0$ and  $l > 0$
the Diophantine system
\begin{eqnarray*}
(n-1)\,i + n\,j & = & k \\
i+j & = & l \\
n & > & 0 \\
i & > & 0 \\
j & \geq & 0
\end{eqnarray*}
has exactly one solution $(n,\,i,\,j) = (\left\lfloor\frac{k}{l}\right\rfloor+1,\,\left\lfloor\frac{k}{l}\right\rfloor l-k+l,\,
k-\left\lfloor\frac{k}{l}\right\rfloor l)$.
\qed

\bigskip
\noindent
If $\chi(x_m) \in \Z[t,u,v]$ is calculated via the recursion of Theorem
\ref{theorek} and if the monomials $t^iu^jv^k$ are ordered by the lexicographical
order of the pairs $(-i,-k)$, we make an important observation: for all $n \geq 1$
the leading term of $\chi(x_{2n+1})$ is $-3\cdot 2^{n-1}\cdot\beta(x_{2n+1})
= -3\cdot (2u)^{n-1}v$.
This assertion is true for $\chi(x_3) = -3v + 3tu + t^3$ and we obtain
for $n \geq 1$ per induction:
$$\chi(x_{2n+3}) = 2u\,\chi(x_{2n+1}) + 2v\,\chi(x_{2n}) + t\cdot(\ldots)
= -3\cdot(2u)^nv + v^2\cdot(\ldots) + t\cdot(\ldots) $$
The ordering of monomials is multiplicative, so (up to a non-zero scalar)
the leading term of $\chi(x_{2n+1}^ix_{2n+3}^j)$ is
$\beta(x_{2n+1}^ix_{2n+3}^j)$.

\begin{lemma}
\label{lemindep}
The monomials $\{\; 1,\, t^i,\, t^i\,x_{2n+1}^{},\, t^i\,x_3^{}\,x_{2n+1}^{},\, x_{2n+1}^i,\,x_{2n+1}^i\,x_{2n+3}^j \,\,
\vert\,\, i,j,n > 0\;\}$ are linearly independent in $\Lambda$.
\end{lemma}
\noindent
{\bf Proof }
Let $\chi_D^{}$ denote the character to $\Q[\sigma_2,\sigma_3]$ corresponding to the Lie
super-algebra $D(2,\!1,\!\alpha)$.
Due to results of \cite{Vo}, one can calculate $\chi_D^{}(x_n)$ with the recursion
of Theorem \ref{theorek} by setting $t = 0$, $u = -2\sigma_2$ and $v = 4\sigma_3$.
The upper observation and Lemma \ref{lembij} make clear that
$\chi_D^{}$ maps the monomials of type $x_{2n+1}^i$ and $x_{2n+1}^i\,x_{2n+3}^j$
to linearly independent elements of $\Q[\sigma_2,\sigma_3]$.

All monomials that contain $t$ are mapped to $0$ by $\chi_D^{}$, so it only remains to
show the linear independence of the monomials of type $t^i,\, t^i\,x_{2n+1},\, t^i\,x_3^{}\,x_{2n+1}$
of a given degree $d$.
We will use the $t = 1$ specialization of the character $\chi_{osp}$ of \cite{Vo}
which we will call $\chi_{osp}$ as well. The target of $\chi_{osp}$ is
$\Q[\alpha]$ and $\chi_{osp}(x_n)$ is given by setting $t = 1$, $u = -\alpha+6\alpha^2$ and
$v = -4\alpha^2+8\alpha^3$ in $\chi(x_n)$.
It is an immediate consequence of the recursion formula that the degree
of the polynomial $\chi_{osp}(x_n)$ is $n$ for $n = 3$ or $n \geq 5$.
So the $\chi_{osp}$-images of $t^d,\, t^{d-2n-1}\,x_{2n+1},\, t^{d-2n-4}\,x_3^{}\,x_{2n+1}$
all have different degrees; this completes the proof.
\qed

\bigskip
\noindent
{\bf Proof of Corollary \ref{corspan}}
Let us introduce an ordering on the monomials $t^{e_0}x_3^{e_1}\ldots x_{2n+1}^{e_n}$
by the lexicographical order on the tuples $(-\sum_{i=0}^n\limits e_i,\,e_1,\,e_2,\,\ldots,\,e_n)$.
The $P$- and $Q$-relations are essentially of the form
\begin{eqnarray*}
3\,x_{i+2}^{}\,x_{j+4}^{} &=& 3\,x_{i+4}^{}\,x_{j+2}^{} + \mbox{ terms with exponent sum} > 2 \\
3\cdot(-2)^{-i}\,x_{2i+2}^{} & = &  \mbox{ terms with exponent sum} > 1 \mbox{\quad(for }i \geq 1).
\end{eqnarray*}
Let $m$ denote a monomial that contains $x_n$ and $x_m$ with $3 \leq n \leq m-3$.
Via the relation $P_{n-2,m-4} = P_{m-4,n-2}$ one can express $m$ as
sum of $m/(x_n\cdot x_m)\cdot x_{n+2}\cdot x_{m-2}$ and terms having a bigger exponent sum.
All monomials in this sum are smaller than $m$.

The $Q$-relations enable us to write any monomial containing $x_{2i}$ in terms
of smaller monomials. Doing both substitutions repeatedly, we obtain 
(after a finite number of steps because the degree is preserved) a linear combination of monomials 
that have only odd indices, which differ at most by $2$
inside each monomial.
This shows that $M$ spans $\Lambda_0$.

To prove the second statement,
let us assume that there is a non-empty linear combination $c$ of
elements of $M$ that is trivial in $\Lambda_0$. Without loss of generality
we may further require that $c$ is homogeneous of degree $d$ with the smallest
possible $d$.
Write $c = t\,\tilde c + r$ where $r$ is a combination of elements of $M$ that do not contain $t$.
Because of $\chi_D^{}(c) = 0$ and $\chi_D^{}(t) = 0$ we have $\chi_D^{}(r) = 0$ and Lemma
\ref{lemindep} implies that $r$ is the empty linear combination.
So $\tilde c$ is a non-empty linear combination of degree $d-1$ and since we assumed
Conjecture \ref{conj2}, $\tilde c$ must be trivial in $\Lambda_0$ contradicting the
minimality of $d$. \qed

\bigskip
\noindent
{\bf Proof of Corollary \ref{coro2} }
For the lower bound we have to estimate the number of monomials of degree
$d$ that appear in Lemma \ref{lemindep}. The number of monomials of the form $t^d$, $t^{d-2n-1}\,x_{2n+1}$ or
$t^{d-2n-4}\,x_3\,x_{2n+1}$ is at least $1 + \left\lfloor\frac{d-2}2\right\rfloor +
\left\lfloor\frac{d-5}2\right\rfloor = d-3$ for $d \geq 0$ (for $d \geq 5$ 
this is even the exact value). 
To complete the proof of the lower bound, 
in view of Lemma \ref{lembij} it remains to show
$$\# \{ \; (i,j) \; \vert \; 2i+3j = d,\; i \geq 0,\; j > 0\; \} \;\geq\; \left\lfloor \frac{d}6 \right\rfloor.$$
\noindent
For any $d$ one may write $d = 6\cdot\left\lfloor\frac{d}6\right\rfloor + 2\,p+3\,q$
with $(p,q) \in \{(0,0),(2,-1),(1,0),(0,1),(2,0),(1,1)\}$.
Then for $1 \leq k \leq \left\lfloor\frac{d}6\right\rfloor$ the pair $(p+3(\left\lfloor\frac{d}6\right\rfloor - k),
\,q+2k)$ is an element of this set.

\medskip
\noindent
For the upper bound, we have to count the monomials of degree $d$ in the set $M$ of Corollary
\ref{corspan}. Lemma \ref{lembij} implies that this number is $1+N(d)$, where
$$N(d) \;\;:=\;\; \#\{\; u^iv^j \; \vert \; 2i+3j \leq d, \; i\geq 0, \; j > 0\;\}.$$
With $t := \left\lfloor\frac{d}{3}\right\rfloor$ and $\varepsilon(n) := \frac{(-1)^n-1}2$ we get:
\begin{eqnarray*}
N(d) & = & \sum_{j=1}^t  \#\{\; i \; \vert \; 0 \leq 2i \leq d-3j \;\}
= \sum_{j=1}^t \left\lfloor\frac{d-3j}{2}+1\right\rfloor
= \sum_{j=1}^t \frac{d-3j+2 + \varepsilon(d+j)}{2} \\
& = & \frac12t\,d-\frac34t(t+1)+t+\frac12\sum_{j=1}^t\varepsilon(d+j) =
\frac12t\,d -\frac34t^2 + \frac14\mu(d)
\end{eqnarray*}
In the last equality we made the abbreviation $\mu(d) := \left\lfloor\frac{d}{3}\right\rfloor+2\cdot
\sum_{j=1}^{\left\lfloor\frac{d}{3}\right\rfloor}\limits \varepsilon(d+j)$.\\
It turns out that
$\mu(d) = 0,0,0,1,-1,1$ when $d \equiv 0,1,2,3,4,5$ modulo $6$.
If we set $r := d - 3\,t$, we finally obtain
$\frac{d^2}{12} + \frac12 - N(d) = \frac{r^2}{12}-\frac{\mu(d)}4 + \frac12$,
which is a value between $0$ and $1$ for all $d$. This shows that $N(d) =
\left\lfloor\frac{d^2}{12} + \frac12\right\rfloor$.\qed
\bigskip

\noindent
{\bf Proof of Theorem \ref{theorek}}
The $\Leftarrow$-direction is easy; we have already seen that
if a $\chi$ with demanded properties exists, it maps the elements
of the spanning set $M$ of type $t^ix_{2n+1}^jx_{2n+3}^k$ to polynomials
with leading terms $t^i\beta(x_{2n+1}^jx_{2n+3}^k)$, which are linearly independent.
Thus $\chi$ is injective and $\Lambda_0$ has no zero divisors because
$\Q[t,u,v]$ is an integral domain.

Let $\hat x_n$ denote the polynomial in $t,u,v$ that is assigned to $x_n$
by the recursive rule in Theorem \ref{theorek}. To prove the $\Rightarrow$-direction
of the statement, we have to show that the $\hat x_n$ satisfy the $P$- and $Q$-relations
(i.e.~replacing all $x_n$ by $\hat x_n$ turns $P$- and $Q$-relations 
into identities in $\Q[t,u,v]$).
For the $P$-relations this can be shown in a straightforward way:
Let $R_{i,j} := P_{i,j}-P_{j,i}$ and $S_i := t\,x_{i-1} + u\,(x_{i-2} - t^{n+1})
+ v\,\big(2\,x_{n} - t^n - 2\,(2t)^n\big)$. With a little patience, one
verifies that replacing every $x_{j+p}$ in $R_{i,j}$ by $S_{j+p}$
(i.e. for $p = 0,1,2,3,4$) leads
to a combination of $R_{i,k}$-s with $k < j$:
$$R_{i,j}\vert_{x_{j} = S_j,\,\ldots\,,\, x_{j+4}^{} = S_{j+4}} \;\;=\;\;
t\;R_{i,j-1} + 2\,u\,R_{i,j-2} + 2\,v\,R_{i,j-3}.$$
So if the $\hat x_n$ satisfy the relation $P_{i,j}-P_{j,i}$ for $j \leq n$, $n \geq 4$,
then $P_{i,n+1}-P_{n+1,i}$ is also satisfied.
To complete the inductive argument, only the cases $(i,j) = (1,2),\, (1,3),\, (2,3)$ have to be checked
by explicit calculation.

\newcommand{\q}[2]{{\hat q}_{#1,#2}^k}
\newcommand{\qk}[3]{{\hat q}_{#1,#2}^#3}
The $Q$-relations are more complicated to handle. Let $\q ij$
denote the polynomial of $\Q[t,u,v]$ that is obtained when in $Q_{i,j}^k$
every $x_n$ is replaced by $\hat x_n$. By definition we have for $i \geq 1$, $j\geq 0$
\begin{equation}
\label{rek1}
\q ij \;=\; t\,\q{i-1}{j+1} - \frac12\q{i-1}{j+2} + \frac12 t^j\q{i-1}2
+ \frac12(\hat x_{j+2}-t^{j+2})\,\q{i-2}2.
\end{equation}
For $i \geq 0$, $j \geq 3$ and $k \in \{1,2\}$ there is another helpful relation, namely
\begin{equation}
\label{rek2}
\q ij \;=\; t\,\q i{j-1} + 2\,u\,\q i{j-2} + 2\,v\,\q i{j-3} - (u\,t^{j-2} + v\,t^{j-3})\,\q{i-1}2
- v\,(2t)^{j-1}\q{i-2}2.
\end{equation}
To make it correct for $i = 0$, we set $\q{-2}2 := 2\,(2t)^{k-2}$. With $\q{-1}2 = t^k$ and
$\q0j = \hat x_{j+k}$ this relation translates into the recursive definition
of the $\hat x_n$ for $i = 0$. To prove the relation for $i \geq 0$, one uses (\ref{rek1}) to
reduce it to terms with first index $< i$ and concludes by induction;
we leave the detailed calculation as an exercise.

The equation (\ref{rek1}) is kindly simple for $j = 0$:
$\,\q i0 \; = \; t\,\q{i-1}1$. Together with (\ref{rek2}) we can thus write
any $\q ij$ in terms of $\q ab$-s with $a \leq i$ and $b \in \{1,2\}$.
When we do this to the right side of equation (\ref{rek1}) for the cases $j=1$ and $j=2$, 
we obtain (for $i \geq 1$):
\begin{eqnarray}
\label{rek3}
\q i1 & = & t\,\q{i-1}2 - u\,\q{i-1}1 + (2\,t\,u-v)\,\q{i-2}2 -  t\,v\,\q{i-2}1
+  2\,t^2v\,\q{i-3}2 \\
\label{rek4}
\q i2 & = & (t^2-u)\,\q{i-1}2 + (t\,u-v)\,\q{i-1}1 + 2\,t\,(t\,u-v)\,\q{i-2}2 -  t^2\,v\,\q{i-2}1
+  2\,t^3v\,\q{i-3}2
\end{eqnarray}
Now we are able to prove that the following two identities
$$
\qk i12 \;\; = \;\; \qk i21, \mbox{\quad and\quad} \qk {i-1}22 \;\; = \;\; \qk i11
$$
hold for $i \geq 0$. It is easy to check them for $i = 0,1,2$.
For $n \geq 3$, one has to assume that {\it both} equations hold for $i < n$
and deduce the case $i = n$ by skillful use of (\ref{rek3}) and (\ref{rek4}).
Once more, this is a lengthy calculation that we will leave to
the reader's ambition.
Up to some nasty computation, we have finally shown that the $\hat x_n$
satisfy the relations $P_{i,j} = P_{j,i}$ and $Q_{i,1}^1 = Q_{i-1,2}^2$.
Under the premise that $t$ is injective, we have already seen that these are
the only relations in terms of $x_n$'s in $\Lambda_0$, so $\chi$ is well defined.

As we have seen above, $\chi$ is then injective, and the recursive formula
for $\chi$ allows us to verify that $\chi(x_n)$ is always of the
form $a + (t^3-tu+v)b$ with $a \in \Q[t]$ and $b \in \Q[t,u,v]$. 
So finally, Corollary \ref{corspan} and Lemma \ref{lembij} imply
the last statement of Theorem \ref{theorek}, which completes the proof.
\qed

\begin{verbatim}
 
e-mail:jan@kneissler.info
http://www.kneissler.info
\end{verbatim}


\begin{thebibliography}{00}

\bibitem{BN} D.~Bar-Natan, {\it On the Vassiliev knot invariants,} Topology {\bf 34} (1995), 423-472.

\bibitem{Kn} J.~A.~Kneissler, {\it The number of primitive Vassiliev invariants
up to degree twelve,} e-print archive (\verb|http://xxx.lanl.gov|), \verb.q-alg/9706022..

\bibitem{cc1} J.~A.~Kneissler, {\it On spaces of connected graphs I: Properties of ladders}, 
Proc.~Internat.~Conf.~"Knots in Hellas '98", Series on Knots and Everything, 
vol.~{\bf 24} (2000), 252-273. 

\bibitem{cc3} J.~A.~Kneissler, {\it On spaces of connected graphs III: The
ladder filtration}, Jour.~of Knot Theory and its Ramif.~Vol.~{\bf 10}, No. 5 (2001), 675-686.



\bibitem{Va} V.~A. Vassiliev, {\it Cohomology of knot spaces,} Theory of Singularities and its Applications
(ed. V.~I. Arnold), Advances in Soviet Math., {\bf 1} (1990), 23-69.

\bibitem{Vo} Pierre Vogel, {\it Algebraic structures on modules of diagrams,} Universit\'e Paris VII preprint, July 1995 (revised 1997).

\bibitem{Vo2} Pierre Vogel, {\it The universal Lie algebra,} Universit\'e Paris VII preprint, June 1999.

\end{thebibliography}
\end{document}